\documentclass[12pt]{article}

\usepackage{amsmath}
\usepackage{amsthm}
\usepackage[pdftex]{graphicx}
\usepackage[all]{xy}
\usepackage{ccaption}

\newtheorem{theorem}{Theorem}[section]

\newtheorem{lemma}[theorem]{Lemma}

\newtheorem{conjecture}[theorem]{Conjecture}

\theoremstyle{remark}

\theoremstyle{definition}
\newtheorem{definition}[theorem]{Definition}
\newtheorem{example}[theorem]{Example}

\newcommand{\script}[1]{\text{$\cal{#1}$}}
\newcommand{\comment}[1]{}
\newcommand{\greg}[1]{}
\newcommand{\rot}[1]{\text{rot}(#1)}
\newcommand{\bill}[1]{}
\usepackage[margin=1in]{geometry}        

  \captionnamefont{\bfseries}
  \captiontitlefont{\small\sffamily}
  \captiondelim{ --- }
  \hangcaption

\begin{document}

\title{On Universal Cycles of Labeled Graphs}
\author{Greg Brockman, Bill Kay, and Emma E. Snively\\
\small Harvard University, University of South Carolina, Rose-Hulman Institute of Technology
}

\maketitle

\begin{abstract}
A universal cycle is a compact listing of a class of combinatorial objects.  In this paper, we prove the existence of universal cycles of classes of labeled graphs, including simple graphs, trees, graphs with $m$ edges, graphs with loops, graphs with multiple edges (with up to $m$ duplications of each edge), directed graphs, hypergraphs, and $k$-uniform hypergraphs.\greg{redid intro}
\end{abstract}

\section{Introduction}

A simple example of a \textit{universal cycle} (U-cycle) is the cyclic string $11101000$, which contains every 3-letter word on a binary alphabet precisely once.  We obtain these words by taking substrings of length 3; it is useful to imagine that we are looking at the string through a ``window'' of length 3, and we shift the window to transition from one word to the next, allowing the window to wrap if necessary. \greg{edited this}

Universal cycles have been shown to exist for words of any length and for any alphabet size.  (For the special case of a binary alphabet, such strings are also known as \textit{de Bruijn cycles}).  The concept easily lends itself to extension, and universal cycles for permutations, partitions, and certain classes of functions are well-studied in the literature (see Chung, Diaconis, Graham~\cite{chungdiaconisgraham} for an overview of previous work in the field\greg{edited this sentence}).  In all cases, the distinguishing feature of a universal cycle is that by shifting a window through a cyclic string (or in some generalizations, an array), all objects in a given class are represented precisely once.

In this paper we generalize the notion of universal cycles.  In particular, we show that these cycles exist for certain classes of labeled graphs. In order to define a universal cycle of graphs, we must first extend the notion of a ``window.'' 

\begin{definition}
Given a labeled graph $G$ with vertex set $V(G)=\{v_0,v_1,\ldots,v_{n-1}\}$ with vertices labeled by the rule ${v_j \mapsto j}$ and an integer $0\leq k \leq n-1$, define a \textit{$k$-window} of $G$ to be the subgraph of $G$ induced by the vertex set $V = \{v_i, v_{i+1}, \ldots, v_{i+k-1}\}$ for some $i$, where vertex subscripts are reduced modulo $n$ as appropriate, and vertices are relabeled such that $v_i\mapsto 0, v_{i+1}\mapsto 1, \ldots, v_{i+k-1}\mapsto k-1$.  For each value of $i$ such that $0\leq i\leq n-1$, we denote the corresponding \textit{$i^\text{th}$ $k$-window} of $G$ as $W_{G,k}(i)$.  If $G$ is clear from context, we abbreviate our window as $W_k(i)$.
\end{definition}

\begin{figure}[ht]
\centering
\includegraphics{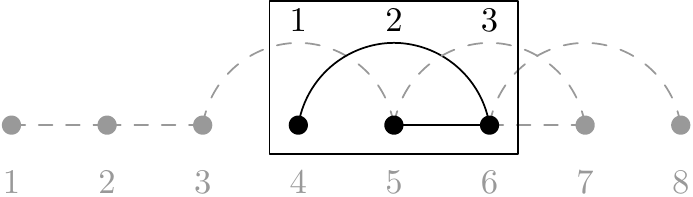}
\caption{A length 3 window of an 8 vertex graph.}
\label{fig:kwindow}
\end{figure}

\begin{definition}
Given $\script F$, a family of labeled graphs on $k$ vertices, a \textit{graph universal cycle (GU-cycle) of $\script F$}, is a labeled graph $G$ such that the sequence of $k$-windows of $G$ contains each graph in $\script F$ precisely once.  That is, $\{W_k(i) | 0 \leq i \leq n-1\} = \script F$, and $W_k(i) = W_k(j) \Rightarrow i = j$.  (Note that the vertex set of the $k$-windows and the elements of $\script{F}$ may be different, however, we will set two labeled graphs equal if they differ only by a bijection between their vertex sets.) 
\end{definition}

\begin{example}
Note that the full 8 vertex graph in Figure \ref{fig:kwindow} is a U-cycle of simple labeled graphs (graphs without loops or multiple edges) on 3 vertices.
\end{example}

\section{Universal cycles of simple labeled graphs}
\label{simple graphs}

  
We begin our investigation by considering only simple graphs; that is, those without loops or multiple edges.  To achieve a form more suited for discourse, we will represent our graphs as strings of integers.  The following definition makes this precise.

\begin{definition}
\label{simple graph encoding}
Consider any simple labeled graph $G$ on $k$ vertices.  The \textit{encoded form} of $G$ is a listing of $k-1$ integers $a_0 a_1 \ldots a_{k-2}$ chosen as follows:

\begin{center}
\begin{tabular}{p{5in}}
For each $0\leq i\leq k-2$, set the base 2 representation of $a_i$ to $b_{k-1} b_{k-2} \ldots b_{i+1}$, where $b_m = 1$ if and only if the vertex with label $i$ is adjacent to the vertex $m$ in $G$.
\end{tabular}
\end{center}

The $a_i$ are the \textit{entries} of $G$'s encoded form.
\end{definition}

\begin{figure}[htp]
\centering
\includegraphics[height=60px]{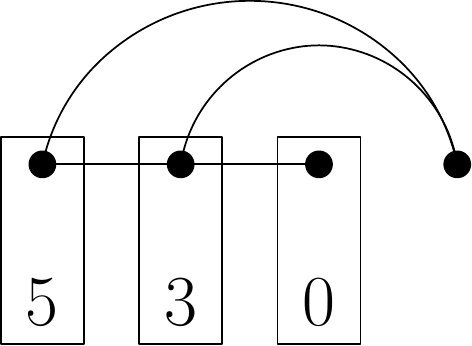}
\caption{A graph and its encoded form, 530}
\label{fig:encoding}
\end{figure}

The following lemma establishes a bijection between simple labeled graphs and their encoded forms.

\begin{lemma}
\label{encodings}
Every simple labeled graph has precisely one encoded form, whose entries $z_i$ satisfy $0 \leq z_i < 2^{k-i-1}$.  Conversely, let $A = a_0 a_1 \ldots a_{k-2}$ be a string of $k$ integers such that $0 \leq a_i < 2^{k-i-1}$.  Then $A$ is the encoded form of precisely one simple labeled graph on $k$ vertices.
\end{lemma}

\begin{proof}
Let $G$ be a simple labeled graph on $k$ vertices with encoded form $a_0 a_1\ldots a_{k-2}$.  By Definition \ref{simple graph encoding}, we see that each $a_i$ has a fixed binary representation.  Furthermore, each $a_i$ is at most ${\underbrace{111\ldots 1}_{k-i-1\text{ 1s}}}_2 = 2^{k-i-1}-1$.  Thus, $0 \leq a_i < 2^{k-i-1}$.

We calculate the total number of strings of integers $a_0 a_1 \ldots a_{k-2}$ such that $0\leq a_i<2^{k-i-1}$ for each $i$ as
$$
\begin{array}{lll}
\displaystyle \prod_{i=0}^{k-2}{2^{k-i-1}} & = & \displaystyle \prod_{i=1}^{k-1} {2^i} \\
                         							     & = & \displaystyle 2^{1+\ldots+(k-1)}\\
                      					    		   & = & \displaystyle 2^{\binom{k}{2}}
\end{array}
$$
This is precisely the number of simple labeled graphs on $k$ vertices.  Since we have established that each such graph has precisely one encoded form, it follows that each such string is the encoded form of exactly one graph.
\end{proof}

Notice that if a graph $G$ on $n$ vertices has encoding $a_0 a_1 \ldots a_{n-1}$, then any $k$-vertex subgraph having consecutive vertices labeled $i, i+1, \ldots, i+k-1$ in $G$ will have encoded form of $a_i \pmod{2^{k-1}}\ a_{i+1} \pmod{2^{k-2}} \ldots\ a_{i+k-2} \pmod{2^1}$ (where subscripts and labels are reduced modulo $n$).  Thus, in a sense, the string $a_i a_{i+1} \ldots a_{i+k-1}$ is equivalent to an encoding with each entry reduced modulo an appropriate power of two; this gives motivation for defining a U-cycle of encodings.

\begin{definition}
\label{equivalence}
Given two strings of $k$ integers, $A = a_0 a_1 \ldots a_{k-1}$ and $B = b_0 b_1 \ldots b_{k-1}$, we say that $A$ and $B$ are \textit{encoding equivalent} if $a_i \equiv b_i \pmod{2^{k-i-1}}$ for each integer $0\leq i\leq k-1$.  Clearly encoding equivalence is an equivalence relation.  Let the canonical representative of each equivalence class be the member $C = c_0 c_1 \ldots c_{k-1}$ with $0 \leq c_i< 2^{k-i-1}$.  Notice that $C$ is the only member of $[C]$ which is itself the encoding of a graph.  If $A$ is a string, we let $C(A)$ denote the graph encoded by the canonical representative from $A$'s equivalence class.
\end{definition}

\begin{definition}
Given a string of integers $A = a_0 a_1 \ldots a_{n-1}$ and a positive integer $k \leq n-1$, define a \textit{$k$-window} of $A$ to be a length-$k$ substring of the form $n_i n_{i+1} \ldots n_{i+k-1}$, where subscripts are reduced modulo $k$ as appropriate.  A \textit{U-cycle of encodings of length $k$} is a string of positive integers such that the sequence of $k$-windows of the string contains exactly one member from each equivalence class of encodings of length $k$.
\end{definition}

\begin{lemma}
A U-cycle of encodings of length $k$ exists if and only if there exists a U-cycle of simple labeled graphs on $k+1$ vertices.
\label{Graphifencode}
\end{lemma}

\begin{proof}
Suppose that $N = n_0 n_1 \ldots n_{l-1}$ is a U-cycle of encodings of length $k$.  Then let $G$ be the graph $C(N)$. Notice that $i^\text{th}$ $k$-window of $N$ encodes the $i^\text{th}$ $k+1$-window of $G$. Since the collection of $k$-windows of $N$ contains precisely one member from each equivalence class of encodings, we see that the collection of $k+1$-windows of $G$ will contain each possible graph on $k+1$ vertices exactly once.

Similarly, let $G$ be a U-cycle of graphs on $k+1$ vertices.  Let $N = n_0 n_1 \ldots n_{l-1}$ be its encoding.  Since taking every $k+1$-window of $G$ yields every $k+1$ vertex graph precisely once, we see that taking every $k$-window of $N$ will yield exactly one element from each equivalence class of encodings.
\end{proof}

We now prove that U-cycles of simple labeled graphs exist for all $k \geq 0$, $k\neq 2$.  We employ two common notions from the study of U-cycles: the transition graph and arc digraph.  The transition graph $T$ of a family $\script{F}$ of combinatorial objects is a directed graph with vertex set $V(T) = \script{F}$. Two vertices $A,B$ are adjacent in $T$ if and only if $B$ can follow $A$ in one window shift of a U-cycle.  If $\script{F}$ is a family of graphs on $k$ vertices, this means that the subgraph induced by the nodes labeled $1, 2, \ldots, k-1$ in $A$ is equal to that induced by the nodes $0, 1, \ldots, k-2$ in $B$.  It should be clear that a U-cycle of $\script F$ corresponds to a Hamiltonian circuit in $T$ (a directed path passing through every node exactly once).

Unfortunately, finding Hamiltonian circuits in graphs is an NP-hard problem; however, in our case the problem can be further reduced.  Let $D$ be the graph with $E(D) = \script{F}$ such that two edges $A,B$ in $D$ are incident if and only if $B$ can follow $A$ in a U-cycle (note that $V(D)$ is arbitrary).  Call $D$ the \textit{arc digraph} of \script{F}.  Now finding a U-cycle of \script{F} is equivalent to finding an Eulerian circuit in $D$ (a directed cycle passing through every edge exactly once); such circuits are easy to detect.  For convenience, we often choose the vertices in $V(D)$ to be labeled with the ``overlap'' between adjacent edges.

\begin{lemma}
\label{simple graph de bruijn}
The arc digraph $D$ of encodings of graphs on $k+1$ vertices, $k\geq 2$, has the following properties:

\begin{enumerate}
\item \label{degrees equal} For each $\script X\in V(D)$, the in-degree of $\script X$ equals the out-degree of $\script X$.
\item \label{strong connectedness} The graph $D$ is strongly connected (there is a directed path from $\script X$ to $\script Y$ for any $\script X\neq\script Y\in V(D)$).
\end{enumerate}
\end{lemma}

\begin{proof}
Fix $k$.  We begin by constructing the transition graph of equivalence classes of encodings of length $k$.  The vertices of this graph are our equivalence classes.  (We provide an example of such a graph in Figure \ref{fig:transition graph}.)

\begin{figure}
\centering
\begin{tabular}{lp{1in}r}
\includegraphics[height=200px]{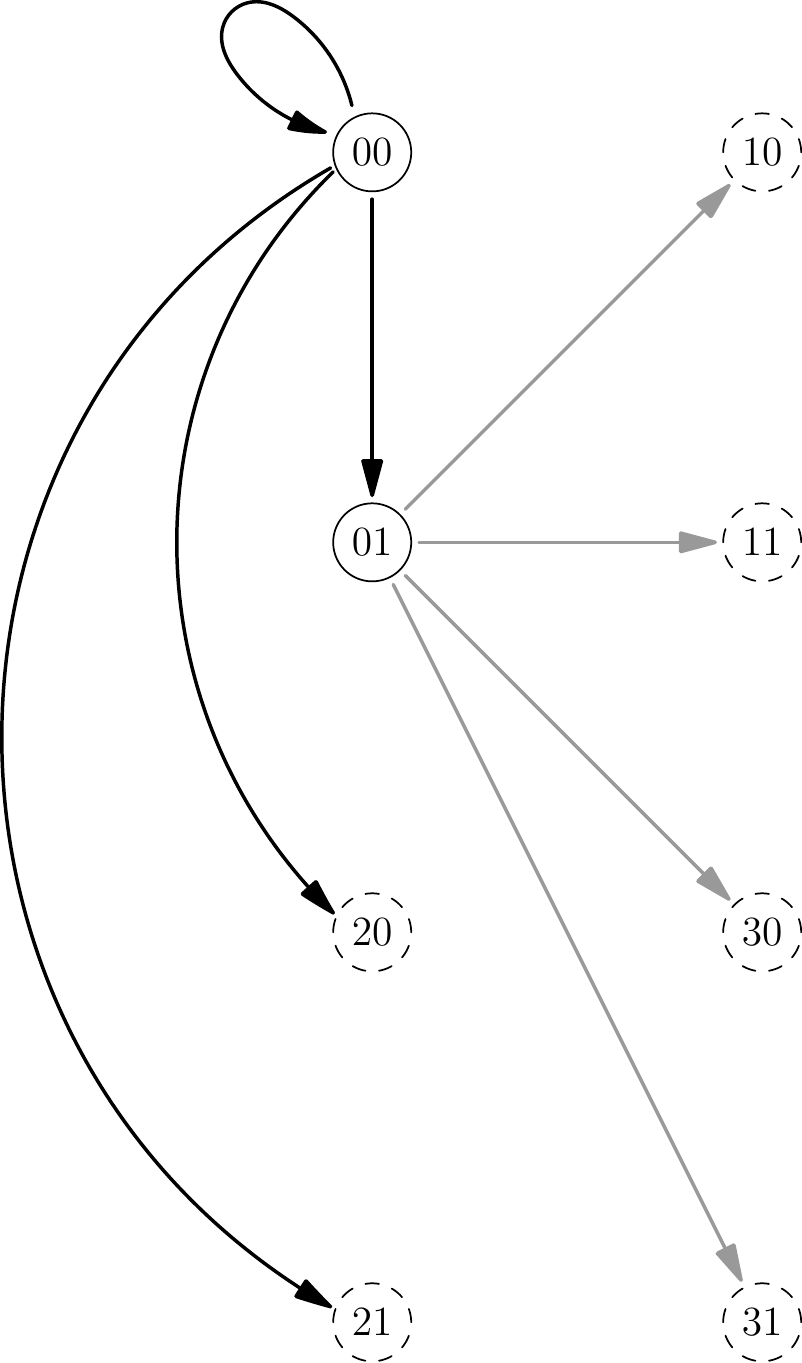}
&
&
\includegraphics[height=200px]{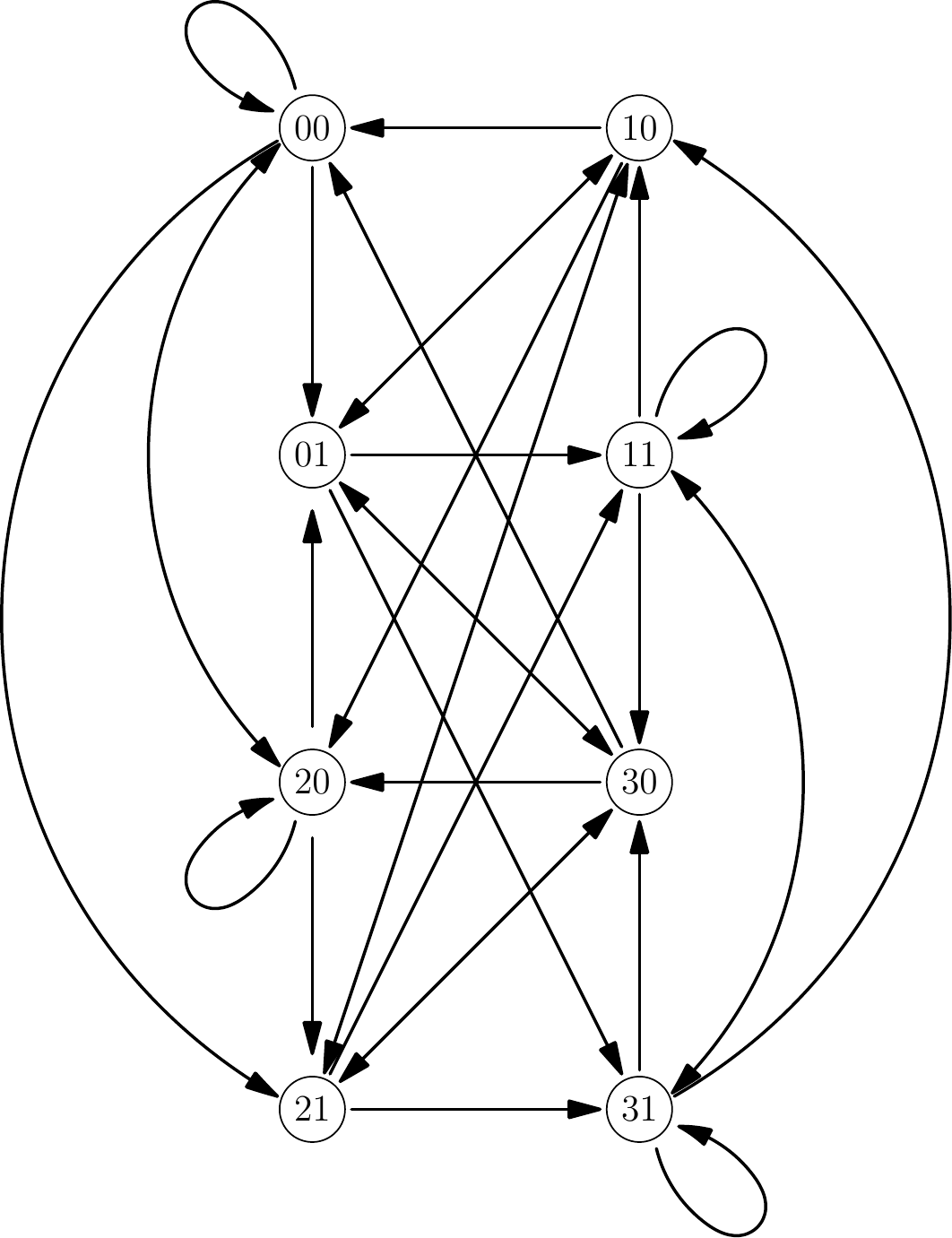}
\end{tabular}
\caption{(Left) A partial sketch of the transition graph of encodings of simple graphs on 3 vertices, and (right) the full transition graph.  We provide the left figure for clarity.}
\label{fig:transition graph}
\end{figure}

Given two encodings with representatives $A = a_0 a_1 \ldots a_{k-1}, B = b_0 b_1 \ldots b_{k-1}$, we draw an edge from $A$ to $B$ if and only if $B$ can follow $A$ in a U-cycle.  This occurs precisely when we can find a member $A' = a_0' a_1' \ldots a_{k-1}'$ of $A$'s equivalence class and a member $B' = b_0' b_1' \ldots b_{k-1}'$ of $B$'s equivalence class such that $a_1' = b_0', \ldots, a_{k-1}' = b_{k-2}'$.  This is equivalent to 
$a_i \equiv b_{i-1} \pmod{2^{k-i}}$ for $1\leq i\leq k-1$.

Consider now the arc digraph corresponding to this transition graph.  Its edge set will be the vertex set of our transition graph.  We will choose the label of each edge to be the canonical representative of that edge's equivalence class.  Thus, we may label our graph's vertices as strings of the form $\alpha_0 \alpha_1 \ldots \alpha_{k-2}$, where $0 \leq \alpha_i < 2^{k-i-1}$.  (Alternatively, we may think about each $\alpha_i$ as being a residue modulo $2^{k-i-1}$.)  See Figure \ref{fig:arc digraph} for an example.

\begin{figure}
\centering
\includegraphics[height=100px]{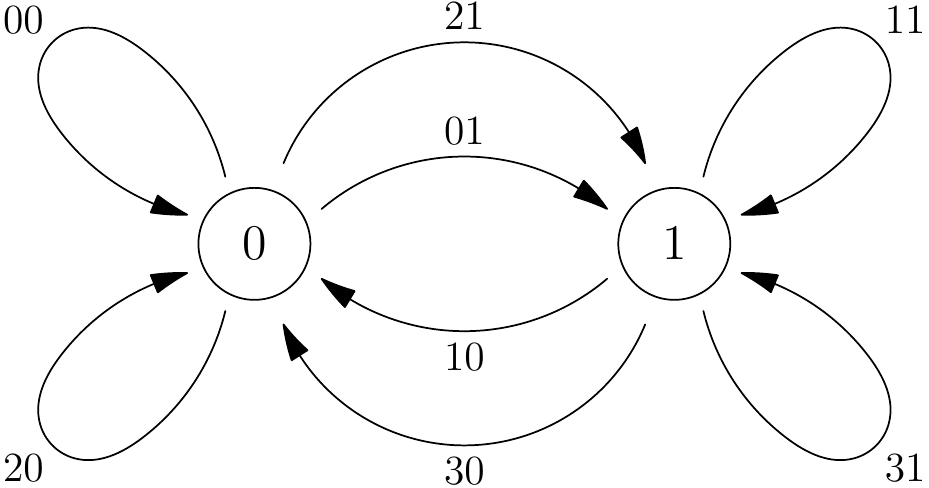}
\caption{The arc digraph for equivalence classes of encodings of length 2 (equivalently, the arc digraph of simple labeled graphs on 3 vertices).}
\label{fig:arc digraph}
\end{figure}

Let $\script{X}, \script{Y} \in V(D)$. Suppose that $A = a_0 a_1 \ldots a_{k-1}$ is an edge of this arc digraph that points from $\cal{X}$ to $\cal{Y}$.  Let $x_i \equiv a_i \pmod{2^{k-i-1}}$ such that $0 \leq x_i < 2^{k-i-1}$; we then have that $\script{X} = x_0 x_1 \ldots x_{k-2}$ and $\script{Y} = x_1 x_2 \ldots x_{k-1}$.

\textbf{Proof of \ref{degrees equal}:} Consider any $\script X\in V(D)$.  Let $\script X = x_0 x_1 \ldots x_{k-2}$.  Then we see that an edge $A = a_0 a_1 \ldots a_{k-1}$ has tail $\script X$ if and only if $x_i \equiv a_i \pmod{2^{k-i-1}}$ for each $0\leq i\leq k-2$.  Now for each $0 \leq i \leq k-2$, we have that $0\leq a_i < 2^{k-i}$ and $0 \leq x_i< 2^{k-i-1}$.  Thus the possible values for $a_i$ are $x_i, x_i + 2^{k-i-1}$.  Hence there are 2 choices for each $a_i$, and there are $2^{k-1}$ possible edges $A$.  Thus $\script X$ has out-degree $2^{k-1}$.

Similarly, an edge $A = a_0 a_1 \ldots a_{k-1}$ has head $\script X$ if and only if $x_i \equiv a_{i+1} \pmod{2^{k-i-1}}$ for each $0 \leq i \leq k-2$.  For each $0\leq i \leq k-2$, we then have that $0 \leq x_i < 2^{k-i-1}$, and it follows that there is only one choice for $a_{i+1}$.  However, $a_0$ can now be chosen arbitrarily, so long as $0\leq a_0<2^{k-1}$.  Hence there are $2^{k-1}$ possible edges $A$, and $\script X$ has in-degree $2^{k-1}$, as desired.

\textbf{Proof of \ref{strong connectedness}:} Consider any two vertices of $D$, $\script X = x_0 x_1 \ldots x_{k-2}$ and $\script Y = y_0 y_1 \ldots y_{k-2}$.  Define $y_{i,j}$ as the residue of $y_i \pmod{j}$, reduced to the interval $[0,j)$.  We show that the following sequence of vertices defines a path in $D$ from $\script X$ to $\script Y$:

$$\begin{array}{lll}
\script X_0 & = & x_0 x_1 \ldots x_{k-2} = \script X \\
\script X_1 & = & x_1 \ldots x_{k-2} y_{0,2^1} \\
\script X_2 & = & x_2 \ldots x_{k-2} y_{0,2^2} y_{1,2^1} \\
&\vdots\\
\script X_i & = & x_i \ldots x_{k-2} y_{0,2^{i}} \ldots y_{i-1,2^1} \\
&\vdots\\
\script X_{k-2} & = & x_{k-2} y_{0,2^{k-2}} \ldots y_{k-3,2^1} \\
\script X_{k-1} & = & y_{0,2^{k-1}} \ldots y_{k-2,2^1} = \script Y
\end{array}$$

Note that each $\script X_i = x_{i}\ldots x_{k-2} y_{0,2^i} \ldots y_{i-1,2^1}$ is a valid vertex label, since the size of its entries is appropriately bounded.  Also, we can easily check that for each $0 \leq i\leq k-2$, there is an edge from $\script X_i$ to $\script X_{i+1}$.
\end{proof}

\begin{theorem}
For each $k \geq 0$, $k\neq 2$, there exists a universal cycle of simple labeled graphs on $k$ nodes.
\end{theorem}

\begin{proof}
When $k=0$ or $k=1$ the result is trivial.  For $k\geq 3$, Lemma \ref{simple graph de bruijn} implies that the arc digraph of equivalence classes of encodings of length $k-1$ has an Eulerian cycle, and consequently there exists a U-cycle of equivalence classes of encodings.  Then applying Lemma \ref{Graphifencode}, we are done.
\end{proof}  

\begin{figure}
\centering
\includegraphics[height=35px]{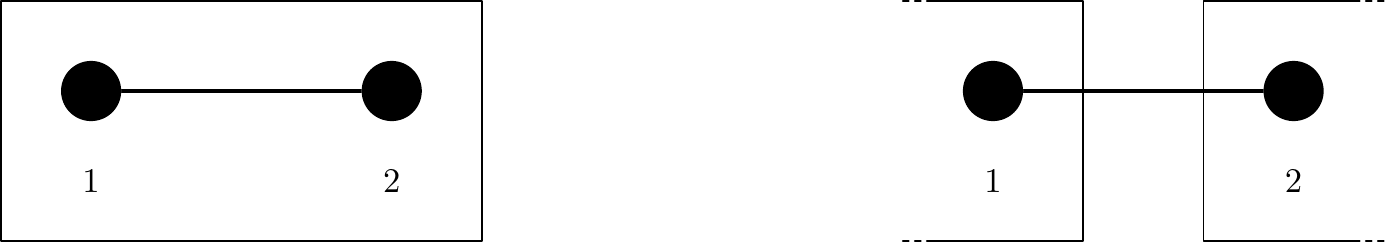}
\caption{An illustration of a GU-cycle using a modified window for $k = 2$.  The left window is the complete graph while the right window is the empty graph.  Note that an edge is considered to be in a window only if it is not ``cut'' by the window.}
\label{fig:k=2}
\end{figure}

Note that for $k=2$ we can modify our definition of a window in order to recognize 2 distinct windows on 2 vertices, as shown in Figure \ref{fig:k=2}.  As a final note to the section, notice that in addition to showing existence, our results provide quick algorithms for constructing the relevant GU-cycles.

\section{General Strategies}

Our approach to showing the existence of GU-cycles in Section \ref{simple graphs} was convenient for rigor; however, we feel that working with encodings is cumbersome.  In the remainder of this text, we will thus avoid the use of encodings in favor of pictorial manipulations.  It should be noted, however, that encodings can be adapted to the situations we are going to consider.

Throughout this section, we suppose that all graphs in a given family have $k$ vertices for some fixed $k$.  Since our results will equally well apply to hypergraphs, we will consider hypergraphs to be a class of graphs.\greg{added this note}

\begin{definition}
Let $\rot X$ be the \textit{rotation class of $X$}, or the set of labeled graphs that differ from $X$ only by a cyclic rotation of vertex labels.
\end{definition}

\begin{lemma}
\label{in equals out}
Let $\script F$ be a family of labeled graphs (including non-simple graphs or even hypergraphs) such that if $X\in\script F$, then $\rot X\subseteq \script F$.  Then in the arc digraph of $\script F$, for every vertex $V$, the in-degree of $V$ equals the out-degree of $V$.\greg{cut this: we have that}
\end{lemma}

\begin{proof}
Let $V$ be a vertex of the arc digraph of $\script F$.  Let the set of edges coming in to $V$ be denoted by $I(V)$, and let the set of edges leaving $V$ be denoted $O(V)$.  We provide a bijection $f: I(V) \longrightarrow O(V)$, thus proving our lemma.  Let $I$ be an edge coming into $V$ (recall that edges in our arc digraph are elements of $\script F$).  If $I$ has $k$ nodes, define $f(I)$ as the graph obtained by cyclically relabeling $I$ as follows: $0\mapsto k-1, 1\mapsto 0, 2\mapsto 1,\ldots, k-1\mapsto k-2$.  Then we see that $f(I) \in\script F$, since $f(I)$ is a rotation of $I$, and furthermore $f(I)$ is an edge leaving $V$.

\begin{figure}[htb]
\centering
\includegraphics{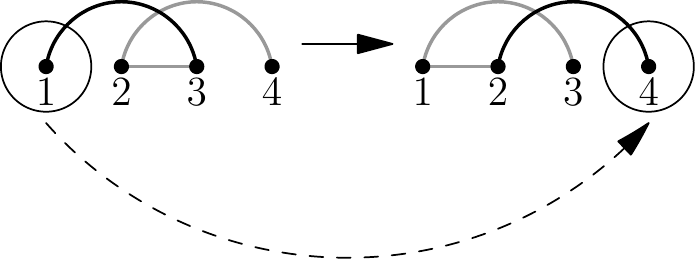}
\caption{The cyclic relabeling of the nodes to create an isomorphic graph}
\label{fig:pushpins}
\end{figure}

Injectivity of $f$ is clear.  Now consider any edge $J$ leaving $V$.  Let $I$ be the graph obtained by cyclically relabeling $J$ as follows: $0 \mapsto 1, 1 \mapsto 2, \ldots, k-2 \mapsto k-1, k-1 \mapsto 0$.  Again, $I$ is isomorphic to $J$, so $I \in\script F$.  Furthermore, $I$ is an edge coming into $V$.  Thus $J$ has a preimage under $f$, and $f$ is surjective.
\end{proof}

Lemma \ref{in equals out} implies that if some class $\script F$ of labeled graphs is closed under rotation, then in order to show that a GU-cycle of \script{F} exists we need only show that the arc digraph is strongly connected (save for isolated nodes). That is, we must show that given two edges $I$ and $J$ in the arc digraph, there exists a directed path in the arc digraph beginning with $I$ and ending with $J$.  In terms of GU-cycles, this is equivalent to showing the existence of a graph $G$ such that $W_{G,k}(i) = I, W_{G,k}(j) = J$ and $W_{G,k}(h) \in \script{F}$ for $i \leq h \leq j$.  Or alternatively, we can picture walking on the arc digraph from $I$ to $J$, taking a series of ``moves'' from one incident edge to another, always following the directed arrows.

We now apply these ideas to prove the existence of U-cycles of various classes of graphs. \greg{ cut this: There exist many broad classes and generalizations of graphs to which we may extend the notion of GU-cycles: non-simple graphs, directed graphs, and hypergraphs. }

\begin{theorem}
\label{extensions of results}
U-cycles exist for the following classes of graphs: graphs with loops, graphs with multiple edges (with up to $m$ duplications of each edge), directed graphs, hypergraphs, and $k$-uniform hypergraphs.
\end{theorem}

\begin{figure}[htb]
\centering
\includegraphics{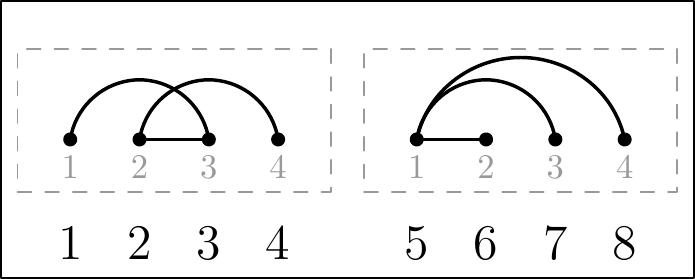}
\caption{Relabeling $J$ to create a new graph $G$}
\label{fig:concat}
\end{figure}

\begin{proof}
Take \script{F} to be any of the above classes of graphs. Pick two graphs $I$ and $J$ from \script{F}. Relabel the nodes of $J$ by the rule $i \mapsto k+i$ to create a graph $G$ such that $V(G) = V(I) \cup V(J)$ and $E(G) = E(I) \cup E(J)$ (see Figure \ref{fig:concat}).  The graph $I$ is the first $k$-window of $G$, and the graph $J$ is the $k^\text{th}$ $k$-window. Further, each $k$-window $W_{G,k}(i)$,  $1 \leq i \leq k$, is a graph in $\script{F}$.  Thus these $k$-windows represent a series of legal edge moves in our arc digraph.
\end{proof}

The extensions from Theorem \ref{extensions of results} followed readily because the relevant graph classes were unrestricted; connectedness of the arc digraph was trivial.  Notice that our proof also applies to some restricted classes of graphs, such as forests.  We now turn our attention to GU-cycles of two types of restricted classes of simple graphs on $k$ vertices.

\begin{theorem}
Graph universal cycles exist for trees on $k$ nodes for $k \geq 3$.
\end{theorem}

\begin{proof}
Let $I$, $J$ be trees.  Let $G$ be the graph obtained by concatenating $I$ and $J$ as in the proof of Theorem \ref{extensions of results}.  As we read the $k$-windows starting from $I$, let $M$ be the first non-tree window that we arrive upon.  Let $v_M$ be vertex of highest label in $M$.  Since none of $G$'s subgraphs contain cycles, we see that there must be one or more components of $M$ that are not connected to the component of $v_M$.  Draw edges from $v_M$ to each of these components; note that the resulting window is now a tree.  Furthermore, these edges did not create any cycles in any of $G$'s $k$-windows, since there are no edges between vertices with label higher than that of $v_M$ and those with lower label.  Also note that still $W_{G,k}(1) = I, W_{G,k}(k) = J$.  We then iterate this process until we arrive at a graph that gives us a sequence of $k$-windows, all of which are trees, starting at $I$ and ending at $J$.
\end{proof}

\begin{theorem}
Graph universal cycles exist for graphs with precisely $m$ edges.
\end{theorem}

\begin{proof}
For any graph $G$, let $d(G)$ be the degree sequence of $G$.  For two graphs $G,H$ having $m$ edges, define $d(G) < d(H)$ if $d(G)$ comes before $d(H)$ lexicographically.  We show that for any graph $I$ having exactly $m$ edges, there is a sequence of moves that takes $I$ to the (unique) graph $L$ having $m$ edges and having least degree sequence of graphs with $m$ edges.  Thereupon, using the bijection we created in Lemma \ref{in equals out}, for any graph $J$ having exactly $m$ edges we can reverse its path to $L$ to arrive at a path from $L$ to $J$.  This will complete our proof.

Let $I$ be a graph with $m$ edges.  If $I$ has least degree sequence, we are done.  Otherwise, let $d(I) = (d_0,d_1,\ldots,d_{k-1})$ and $d(L) = (L_0,L_1,\ldots,L_{k-1})$.  Since $d(I)$ is not minimal, there must be some $i$ such that $d_i > L_i$.  But then there must exist some $j > i$ such that $d_j < L_j$.  Now consider a sequence of moves where at each step we rotate $I$'s vertices according to the relabeling $0\mapsto k-1, 1\mapsto 0,\ldots, k-1\mapsto k-2$ until the vertex formerly labeled $i$ attains label 0; let $i'$ be any vertex adjacent to it.  Once more rotate the vertex set, but this time connect the vertex formerly labeled $i'$ to $j$ instead of $i$.  After rotating all of the vertices back to their original positions, we obtain a graph with a smaller degree sequence than $I$.  By infinite descent, we see that we can eventually move to $L$ via a sequence of moves all of which are graphs having $m$ edges, thus completing the proof of the theorem.
\end{proof}


\section{Conclusions and Future Directions}

In this paper, we have presented a beginning theory of universal cycles of graphs.  We have shown the existence of U-cycles of various classes of labeled graphs on $k$ nodes, including simple graphs, multigraphs, graphs on $m$ edges, directed graphs, trees, hypergraphs, and $k$-uniform hypergraphs.  

\begin{figure}[htb]
\centering
\includegraphics{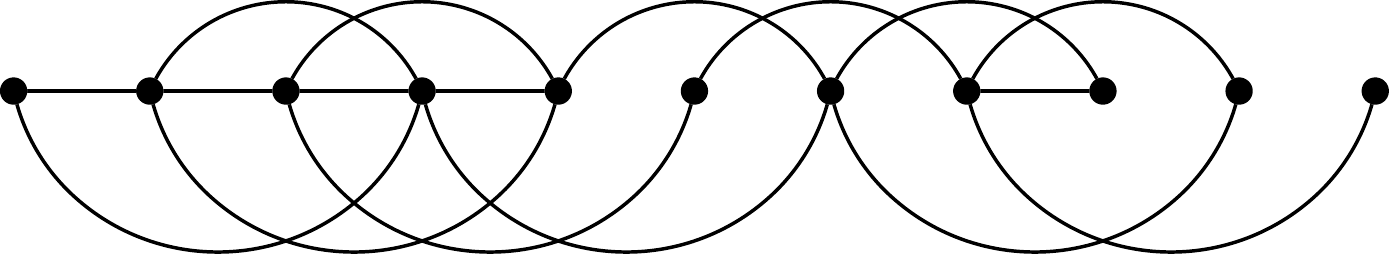}
\caption{A graph universal cycle of isomorphism classes of graphs on 4 nodes.}
\label{fig:unlabeled 4 cycle}
\end{figure}

Our work in this field is far from complete.  There exist many other classes of graphs for which there conceivably exist U-cycles.  However, perhaps the most obvious gap is results regarding U-cycles of unlabeled graphs.  The canonical result would be to prove the existence of U-cycles of isomorphism classes of graphs. In such a U-cycle, no two $k$-windows are isomorphic and every isomorphism class is represented as a $k$-window.  It is easy to find a U-cycle of isomorphism classes of graphs on 3 nodes.  It is difficult, but still possible, to find a U-cycle of isomorphism classes of graphs on 4 nodes; one such cycle in exhibited in Figure \ref{fig:unlabeled 4 cycle}.  These results lead us to conjecture the following.

\begin{conjecture}
For each $k \geq 3$, there exists a U-cycle of isomorphism classes of graphs on $k$ nodes.
\end{conjecture}

We also note that graph universal cycles have potential in theorem-proving as well, as demonstrated by the following result.\greg{edited the wording of this and last paragraph a bit}

\begin{definition}
We say that an integer-valued graph theoretic function $f$ is \textit{window-Lipschitz} if, for all graphs $G$ and $H$ which are one window shift apart in a graph universal cycle, $|f(G) - f(H)| \leq 1$.
\end {definition}

\begin{lemma}
Let $U$ be a U-cycle of some family $\script F$ of graphs, and let $f$ be a window-Lipschitz function defined on these graphs. Then for each integer $\displaystyle\min_{G\in\script  F} f(G) < i < \max_{G\in\script  F} f(G)$ there exist at least two distinct elements of $G\in\script F$ such that $f(G) = i$.
\label{lipschitz}
\end{lemma}

\begin{proof}
By definition of a Lipschitz function, under a single window shift the value of $f$ can change by at most 1.  Hence during the sequence of window shifts from the graph with minimal $f$-value to maximal, every possible value of $f$ in between the minimum and maximum is attained.  Similarly, during the sequence of window shifts from the graph with maximal $f$-value to that with minimal, every possible value of $f$ is again attained.  This completes our proof.
\end{proof}

Some examples of window-Lipschitz functions are chromatic number and largest clique.

\section{Acknowledgements}

This work was done at the East Tennessee State University REU, NSF grant 0552730, under the supervision of Dr. Anant Godbole.

\bibliographystyle{amsplain}
\bibliography{gucyclesbib}

\begin{tabular}{p{3in}l}
\small\textsc{Greg Brockman} & \small\textsc{Bill Kay}\\[-5pt]
\small\textsc{Harvard University} & \small\textsc{University of South Carolina}\\[-5pt]
\small\textsc{Cambridge, MA 02138} & \small\textsc{Columbia, SC 29208}\\[-5pt]
\small\textsc{United States} & \small\textsc{United States}\\[-5pt]
\small\verb|gbrockm@fas.harvard.edu| & \small \verb|kayw@mailbox.sc.edu|
\end{tabular}
\vspace{.3in}

\begin{tabular}{p{3in}l}
\small\textsc{Emma E. Snively}\\[-5pt]
\small\textsc{Rose-Hulman Institute of Technology}\\[-5pt]
\small\textsc{Terre Haute, IN 47803}\\[-5pt]
\small\textsc{United States}\\[-5pt]
\small\verb|snivelee@rose-hulman.edu|
\end{tabular}\\

\end{document}